\documentclass[12pt]{extarticle}
\pdfoutput=1

\usepackage[ngerman,english]{babel}
\usepackage{graphicx}
\usepackage{framed}
\usepackage[normalem]{ulem}
\usepackage{xfrac}
\usepackage{enumerate}
\usepackage{comment}
\usepackage[utf8]{inputenc}
\usepackage[top=1 in,bottom=1in, left=1 in, right=1 in]{geometry}

\usepackage{amsmath}
\usepackage{amsfonts}
\usepackage{amsthm}
\usepackage{amssymb}
\usepackage{mathtools}
\usepackage{hyperref}
\usepackage[utf8]{inputenc} % allow utf-8 input
\usepackage[T1]{fontenc}    % use 8-bit T1 fonts
\usepackage{enumitem}
\usepackage{cleveref}

\newcommand{\abs}[1]{\left\vert #1 \right\vert}
\newcommand{\norm}[2][]{\ensuremath{\left\lVert#2\right\rVert_{#1}}}
\newcommand{\nnorm}[2][]{\ensuremath{\lVert#2\rVert^{}_{ #1}}}
\DeclareMathOperator{\tr}{tr}
\DeclareMathOperator{\rank}{rank}

\newcommand{\N}{\mathbb{N}}
\newcommand{\R}{\mathbb{R}}
\newcommand{\C}{\mathbb{C}}
\newcommand{\Q}{\mathbb{Q}}
\newcommand{\Z}{\mathbb{Z}}
\DeclareMathOperator*{\argmin}{arg\,min}

\theoremstyle{definition}
\newtheorem{definition}{Definition}

\theoremstyle{plain}
\newtheorem{theorem}[definition]{Theorem}
\newtheorem{corollary}[definition]{Corollary}
\newtheorem{lemma}[definition]{Lemma}
\newtheorem{proposition}[definition]{Proposition}

\theoremstyle{remark}
\newtheorem{remark}[definition]{Remark}

\usepackage[autostyle]{csquotes}
\usepackage{biblatex}
\usepackage{authblk}

\title{Turing meets Moore-Penrose: Computing the Pseudoinverse on Turing Machines}
\author[1,2,3]{Holger Boche}
\author[4]{Adalbert Fono}
\author[4,5,6,7]{Gitta Kutyniok}
\affil[1]{Institute of Theoretical Information Technology, TUM School of Computation, Information and Technology, Technical University of Munich, Germany}
\affil[2]{Munich Quantum Valley (MQV)}
\affil[3]{Cluster of Excellence “Centre for Tactile Internet with Human-in-the-Loop” (CeTI) of Technische Universität Dresden, Germany }
\affil[4]{Department of Mathematics, Ludwig-Maximilians-Universität München, Germany}
\affil[5]{Department of Physics and Technology, University of Tromsø, Norway}
\affil[6]{Munich Center for Machine Learning (MCML), Munich, Germany}
\affil[7]{DRL-German Aerospace Center, Germany}
\date{} 

\begin{document}

\maketitle

\begin{abstract}
The pseudoinverse of a matrix, a generalized notion of the inverse, is of fundamental importance in linear algebra and, thereby, in many different fields. 
Despite its proven existence, an algorithmic approach is typically necessary to obtain the pseudoinverse in practical applications. Therefore, we analyze if and to what degree the pseudoinverse can be computed on perfect digital hardware platforms modeled as Turing machines. For this, we utilize the notion of an effective algorithm that describes a provably correct computation: upon an input of any error parameter, the algorithm provides an approximation within the given error bound with respect to the unknown solution. We prove that a universal effective algorithm for computing the pseudoinverse of any matrix with a finite error bound does not exist on Turing machines. 
However, for specific classes of matrices, we show that provably correct algorithms exist and obtain a characterization of the properties of the input set, leading to the effective computability breakdown.
\end{abstract}

%\textbf{Keywords}: pseudoinverse, effective algorithm, Turing machine 

%\textbf{MSC}: 03D10, 3D80, 68Q04, 15A09

\section{Introduction}
\label{sec:1}

The pseudoinverse of a matrix \cite{Moore20GenInv, Penrose55GenInv} plays an important role in various areas of linear algebra and its applications: the solution and statistical analysis of linear systems \cite{Greville59SolLinSys, Albert1976StatAppPI}, the condition number of a matrix, the orthogonal projection onto the range of a matrix, and matrix decomposition \cite{CAIAFA2010MDecomp} are all connected to the associated pseudoinverse. 
Moreover, computing the condition number of a matrix describing a set of linear equations is a pre-requisite to enable speedups in quantum algorithms estimating features of the solution \cite{Harrow2009QuAlg, Barz2014QuAlg, Cai2013QuAlg}. Going beyond pure linear algebra, the pseudoinverse is also an important tool in areas such as graph theory \cite{BELARDO2014Graph, WANG2017Graph, AZIMI2018Graph} and robotics \cite{hock2017forward, tevatia2000inverse}. 
These examples stress the significance of the pseudoinverse in various fields. In practice, approaches to compute the pseudoinverse comprise strategies ranging from rank decomposition, applying the QR method, singular value decomposition \cite{BenIsrael2003GenInv}, and based on Gaussian elimination \cite{Sheng2010CompMP, Ji2012MPinv}. However, these approaches typically lack theoretically proven error control---in the strict sense formalized in Section \ref{sec:CompAna}---or generality, so universal correctness guarantees are not established despite the proven existence and uniqueness of the pseudoinverse \cite{Moore20GenInv, Penrose55GenInv}. Consequently, the algorithms provide seemingly correct outputs for any input instance and are indeed provably reliable in specific cases, such as well-conditioned inputs \cite{Pan91MPNewton, Soderstrom74MPIter}, but intrinsic failure modes and their severity are unbeknownst to the user. 

Note that closed-form representations of the pseudoinverse are known for universal  \cite{Barata12MPTutorial, BenIsrael66IterAlg, campbell2009generalized} as well as specific classes of matrices \cite{Bajo21MPPoly}. However, it is unclear if they can be implemented on digital hardware with error control, i.e., if the actual execution of the involved mathematical operations on digital hardware is feasible in the described sense. This is a crucial point since it is straightforward to reason that in a more general theoretical computing concept, i.e., noise-free, analog computing, the pseudoinverse based on the algorithm introduced in \cite{Sheng2010CompMP} can be computed with correctness guarantees as shown in \cite{Boche2025InvProbAnal}. However, this model exceeds the capability of digital computations and, consequently, is not an adequate description of real-world digital computing.
A (potential) computing paradigm beyond classical digital hardware is given by quantum computing
\cite{Nielsen2010QuComp} and theoretical quantum algorithms to compute the pseudoinverse have indeed been proposed: Under the assumption that the singular values of the input matrices are restricted to a certain interval not containing zero, a quantum singular value transformation has been described that enables the computation of the pseudoinverse \cite{Gilyen2019QSTV, Martyn2021QuantumSurvey}. Despite these positive signs, the quantum computing paradigm does not resolve the issues posed by implementations on classical digital hardware. First, whether quantum algorithms exist that do not depend on a specific input requirement is unknown. Second, it is unclear if and when the potential benefits of quantum algorithms will translate into practical, real-world advantages since actual implementations of quantum computers are still missing \cite{Bloch2022PracticalQuantumAdv}. Third, the computational power of digital quantum computers does not exceed the power of classical digital computers, their main advantage is the potential reduction in computational complexity and the associated speed-up in computing time \cite{kitaev2002classical}.

Therefore, we focus on digital computing, currently encompassing virtually any general-purpose computing device as well as high-performance computer, and pose the question of whether a universal algorithm (with error control) to compute the pseudoinverse of any matrix exists or under which conditions it can be theoretically established.

\subsection{Problem Statement}
The pseudoinverse $A^\dagger \in \C^{n\times m}$ of a matrix $A \in \C^{m\times n}$ is characterized by 
\begin{equation*}
    A A^\dagger A = A, \qquad A^\dagger A A^\dagger = A^\dagger, \qquad  (A A^\dagger)^H = A A^\dagger, \quad \text{and} \quad (A^\dagger A)^H = A^\dagger A.
\end{equation*}
These equations uniquely define the pseudoinverse, i.e., exactly one matrix $A^\dagger$ satisfies all four conditions, albeit a simple closed-form expression of $A^\dagger$ is known only in specific instances. For example, if $A$ has full rank such that its columns are linearly independent, then it is easy to verify that $A^\dagger = (A^H A)^{-1} A^H$. In the general case, an explicit description that allows for a straightforward effective computation is missing, i.e., an algorithm that takes any matrix (via arbitrarily accurate rational approximations) and any error parameter as input and computes an approximation of the associated pseudoinverse up to the given error. The need for effective computations is accentuated by fundamental problems whose solution depends on the pseudoinverse. The well-known least squares problem asks us to 
\begin{equation*}
    \text{ minimize } \norm[2]{Ax -b} \text{ over } x \in \C^n \quad \text{ for } A \in \C^{m\times n} \text{ and } b \in \C^m.
\end{equation*}
This optimization problem is convex, and the unique minimizer with minimum Euclidean norm $A^\dagger b \in \C^n$ is obtained via the pseudoinverse \cite{Stewart77PseudoInv}.

Our core question is whether or under which conditions the pseudoinverse and fundamental problems associated with the pseudoinverse, such as the least squares problem and the condition number of a matrix, can be effectively approximated. To that end, we consider the following mappings:

\begin{enumerate}[label=\Roman*., itemsep=0.1cm]
    \item \makebox[\linewidth]{$f^\dag: \C^{m\times n} \to \C^{n \times m}, \quad A\mapsto A^\dagger$} 
    \item \makebox[\linewidth]{$f^{\nnorm{\dag}}: \C^{m\times n} \to \R, \quad A\mapsto \nnorm[F]{A^\dagger}$}      
    \item \makebox[\linewidth]{$f_{\text{lsq}}: \C^{m\times n} \times \C^m \to \R, \quad (A,b) \mapsto \min_{x\in\C^n} \norm[2]{Ax-b}= \nnorm[2]{A A^\dag b - b}$}  
    \item  \makebox[\linewidth]{$f_\text{lsq-m}: \C^{m\times n} \times \C^m \to \C^n, \quad (A,b) \mapsto \argmin\limits_{x\in \C^n: x \text{ is minimizer of } \norm[2]{Ax-b}} \norm[2]{x} = A^\dagger b$}  
    \item \makebox[\linewidth]{$f_\text{lsq-n}: \C^{m\times n}\times \C^m \to \R, \quad (A,b)\mapsto \nnorm[2]{f_\text{lsq-m}(A,b)}= \nnorm[2]{A^\dagger b}$}
    \item  \makebox[\linewidth]{(condition number) $\kappa: \C^{m\times n} \to \R, \quad A\mapsto \nnorm[F]{A} \nnorm[F]{A^\dagger}$} 
\end{enumerate}

\subsection{Contributions}
We analyze the existence of effective algorithms on digital computers modeled by Turing machines approximating the functions I.-VI. Our focus is on the pseudoinverse mapping---since it plays a crucial role in the computation of all the considered functions---and ask whether an effective algorithm exists to compute the pseudoinverse of every matrix $A\in \C^{m\times n}$? We can split this question into two parts and study them separately: First, given a fixed matrix $A\in \C^{m\times n}$, does a Turing machine $\text{TM}_A$ computing the pseudoinverse $A^\dag$ (with error control) exist? Second, does a Turing machine $\text{TM}_U$ taking any matrix $A\in \C^{m\times n}$ as input and constructing $\text{TM}_A$ as output exist?
\begin{itemize}
    \item We provide, despite the existence of the pseudoinverse, a negative answer to the latter yet a positive answer to the former question. In particular, we show in Theorem \ref{thm:non-approx} that no algorithm solving the second part with finite error exists: Any algorithm will make arbitrarily large errors for certain input instances. This is contrasted with the observation in Proposition \ref{prop:comp_values} that a Turing machine computing the associated pseudoinverse for a fixed matrix exists. Therefore, the computability breakdown in the universal case arises at the level of finding said Turing machine for a given input matrix. The construction of the sought Turing machine requires specific information about the input matrix, which is not provided in this general setting. Exemplarily, we analyze and identify this breakdown in the special case of iterative algorithms with convergence guarantees in Corollary \ref{cor:Effect}.
    \item The previous findings indicate that by diminishing the generality of the problem, i.e., restricting the input domain of potential algorithms, effective algorithmic approaches to compute the pseudoinverse may exist since appropriately reducing the input domain corresponds to increasing the available information about the remaining input instances. We demonstrate for two classes of matrices---matrices of full rank in Theorem \ref{thm:FullRank} and matrices with elements bounded away from zero in O
    Proposition\ref{prop:EpsSet}---that the observation is valid by establishing the existence of effective algorithms in both cases. Hence, in practical applications, it is important to assess the given problem and derive information about the feasible input instances for potentially realizing effective solvability.
    \item Via Theorem \ref{thm:class} and its implication, we fully characterize the computability breakdown and describe a tool to check whether a problem with a specific input domain may be amenable to effective approaches. This also closely aligns with the known discontinuities of the mapping on the pseudoinverse. Moreover, we wish to highlight that no algorithm effectively decides whether an arbitrary matrix satisfies the desired input conditions generally. In this sense, a verifiably correct 'red flag' functionality does not exist for the above benign input classes.  
\end{itemize}

The described computational barrier does not arise due to numerical errors in real-world implementations but is an inherent feature of Turing machines. Indeed, the condition number for invertible matrices provides a measure of numerical errors for implementations of the inversion operation on actual digital hardware with floating point arithmetic. 
However, computing the inverse of a fixed matrix can be performed effectively, even on instances with ill-conditioned inputs. Thus, a non-computable/non-approximable property establishes an impossibility result for any realization on a digital machine. In contrast, the condition number only provides a heuristic to decide whether a specific input is benign in certain circumstances.

\subsection{Outline}
First, we concisely introduce computability theory in Section \ref{sec:CompAna}. Then, in Section \ref{sec:Results}, our positive and negative results concerning the effective computation of the functions I.-VI. are presented. 
Finally, the \nameref{sec:proofs} covers the proof of our statements.

\section{Digital Computations}\label{sec:CompAna}

The predominant computing devices today are digital computers. Digital hardware, in principle, allows for the exact computation of specific problems, such as discrete ones. For continuous problems, algorithms on digital hardware only provide an approximation of the solution. Thus, a well-defined notion of an approximate solution is required to evaluate the computed approximation and guarantee its correctness. An approach is to ask for {\em effective} algorithms, where algorithms compute an approximate solution and quantify the error to the (generally unknown) exact solution. In particular, algorithms must be able to provide an approximation within any prescribed error bound. 
A theoretical model for digital computers, which captures the logic of digital computations but neglects real-world limitations such as memory constraints, energy consumption, etc, is given by Turing machines \cite{Turing36Entscheidung}. Any algorithm that can be executed by a real-world (digital) computer can, in theory, be simulated by a Turing machine and vice versa.
Therefore, the existence of effective algorithms for digital computers is studied via Turing machines \cite{Myhill71Noncomp, PourEl97WaveEq, Boche20LTI, Elkouss2018MemoryEC, Schaefer2019TuringMS, Boche20SpecFac, Boche20BandlimitedSignals}.

\subsection{Computable Analysis}

This section presents the necessary definitions from computability theory on continuous domains based on Turing machines. For a comprehensive overview of real-valued computability, we refer to \cite{Soare87RecursivelyES, Weihrauch00CompAnal, Pour-El17Computability, AvigadBrattka14CompAnal}, whereas elementary topics of computability theory, such as recursive functions on discrete domains and the formal definition of Turing machines, can be found in any textbook \cite{computabilitybook}. 

The general paradigm describing digital computation on real numbers introduced by Turing himself \cite{Turing36Entscheidung} is based on the transformation of sequences approximating real numbers with arbitrary precision. 

\begin{definition} \label{def:compnum}
We introduce the following notions of computable objects.
\begin{itemize}
    \item A sequence $(r_k)_{k \in \N}$ in $\Q$ is \textit{computable}, if there exist recursive functions $a, b, s : \N \to \N$ such that $b(k) \neq 0$ and
    \begin{equation*} 
        r_k = (-1)^{s(k)} \cdot \frac{a(k)}{b(k)} \qquad \text{ for all } k \in \N. 
    \end{equation*} 
    \item A sequence $(x_k)_{k \in \N}$ in $\R$ \textit{converges effectively} to $x \in \R$, if there exists a recursive function $e: \N \to \N$ such that
    \begin{equation*}
        \abs{r_k - x} \leq 2^{-N} \quad \text{ for all } N \in \N \text{ and all } k \geq e(N).
    \end{equation*}
    \item A real number $x \in \R$ is \textit{computable}, if there exists a computable rational sequence $(r_k)_{k \in \N}$ converging effectively to $x$. We refer to the sequence $(r_k)_{k \in \N}$ as a \textit{representation} for $x$ and denote the set of computable real numbers by $\R_c$.
    \item A sequence $(x_n)_{n\in\N}$ in $\R$ is \emph{computable}, if there exists a computable double-indexed rational sequence $(r_{n,k})_{n, k\in\N}$ and a recursive function $e: \N \times \N \to \N$ such that
    \begin{equation*} 
        \abs{x_{n} - r_{n,k}} \leq 2^{-N} \quad \text{ for all } n,N \in \N \text{ and all } k \geq e(n,N).
    \end{equation*}
\end{itemize}
\end{definition}
    
\begin{remark}
    The previous definitions can be extended to $\R^d$ and $\C^d$---considered as a $2d$-dimensional real vector space---by requiring that each one-dimensional component(-wise sequence) is computable. We denote the set of \textit{computable complex numbers} by $\C_c := \{x+iy\in\C : x,y\in\R_c\}$. 
\end{remark}

The computability of functions is a well-studied property, and various computability notions exist \cite{AvigadBrattka14CompAnal}; we will employ the following ones to establish positive and negative results regarding effective computations on digital hardware.

\begin{definition} 
    For $D \subset \C^d$, a function $f: D \to \C_c^m$ is 
    \begin{itemize}
        \item \textit{Borel-Turing computable}, if there exists a Turing machine that transforms each representation of $x \in D\cap \C_C^d$ into a representation for $f(x)$.        
        \item \textit{Banach-Mazur computable}, if $f$ maps any computable sequences $(x_n)_{n\in\N}$ in $D\cap \C_c^d$ onto computable sequences $(f(x_n))_{n\in\N}$.
    \end{itemize}
\end{definition}

\begin{remark}\label{rm:BTc}
    Borel-Turing computability formalizes an intuitive notion of computation, where an algorithm approximates a given function to any accuracy with error control. Thus, effective algorithms can only exist if the associated problem is Borel-Turing computable. A weaker form of computability is given by Banach-Mazur computable functions, which preserve computability of sequences. It is well known that Borel-Turing computability implies Banach-Mazur computability and that computable functions in both senses are necessarily continuous \cite{AvigadBrattka14CompAnal}.
\end{remark}
\begin{remark}\label{rm:elementaryfct}
    Basic arithmetic operations, as well as elementary operations such as the absolute value function, are computable \cite{Pour-El17Computability}. This also immediately implies the computability of the Euclidean and Frobenius norm.
\end{remark}   
Another slightly different notion of computation is established via oracles providing representations for any (and not just computable) objects \cite{Ko91ComplTheoryRealFunc}. However, for our purposes of analyzing digital computations, the natural domain is indeed characterized by computable objects, and we apply the introduced computability framework in the remainder.

\section{Main Results}\label{sec:Results}
 
Next, we present our results regarding the effective computation of the functions I.-VI. The crucial step therein is the computation of the pseudoinverse, and we first show that an effective approach to approximate the pseudoinverse does not generally exist. Subsequently, we analyze the breakdown of existing algorithmic strategies on the basis of iterative processes and establish conditions on the input domain to enable the introduction of effective algorithms for computing the pseudoinverse. The proofs of our findings are provided in the \nameref{sec:proofs}. 

\subsection{Banach-Mazur Non-Computability of the Pseudoinverse}\label{sec:BMnonComp}
Assuming Borel-Turing computability of the mapping of a matrix on its pseudoinverse immediately implies the effective computation of the functions I.-VI. on digital hardware. Can we find a constructive algorithmic approach that enables us to verify the computability of said mapping? For instance, the singular value decomposition (SVD) is a standard method to compute the pseudoinverse. A crucial step therein consists of identifying the non-zero singular values of the input matrix. However, comparisons to zero can not be performed effectively on digital hardware \cite{Pour-El17Computability}. Hence, SVD does not suffice to establish an effective approach to compute the pseudoinverse. 
The next theorem shows that this is not due to a specific property of SVD, but an effective approach to approximate the pseudoinverse does not exist in full generality.

\begin{theorem}\label{thm:non-approx}
    For any function $f$, denote a function with the same domain and codomain by $\hat{f}$. Then, for $m,n \geq 2$, the functions I.-VI. are not algorithmically approximable, i.e., any function $\hat{h}$ is not Banach-Mazur computable if
    \begin{enumerate}[label=(\roman*)]
        \item 
        \begin{equation*}
            \sup_{\substack{A \in \C^{m\times n}: \\ \norm[F]{A} \leq \sqrt{2}}} \nnorm[F]{h(A) - \hat{h}(A)} < \infty \quad \text{ for } h \in \{f^\dag, f^{\nnorm{\dag}}, \kappa\}; 
        \end{equation*}
        the claim is also valid if any other norm replaces $\norm[F]{\cdot}$ in $f^{\nnorm{\dag}}$ and $\kappa$.
        \item 
        \begin{equation*}
            \sup_{\substack{(A,b) \in \C^{m\times n} \times \C^m : \\ \norm[F]{A} \leq \sqrt{2},\norm[2]{b} \leq \sqrt{2}}} |h(A,b) - \hat{h}(A,b)| < \frac{1}{4} \quad \text{ for } h =  f_{\text{lsq}}.
        \end{equation*}

        \item 
        \begin{equation*}
            \sup_{\substack{(A,b) \in \C^{m\times n} \times \C^m : \\ \norm[F]{A} \leq \sqrt{2},\norm[2]{b} \leq \sqrt{2}}} \nnorm[2]{h(A,b) - \hat{h}(A,b)} < \infty \quad \text{ for } h \in \{f_\text{lsq-m},f_\text{lsq-n}\}.  
        \end{equation*}
    \end{enumerate} 
\end{theorem}

\begin{remark}\label{rm:BTnc}
    Since the results are also valid in Borel-Turing sense, effective algorithmic computations on digital hardware of the problems posed by I.-VI. in the sense of Remark \ref{rm:BTc} are not feasible for accuracies below the introduced lower bounds. In particular,  except for $f_{\text{lsq}}$, there indeed does not exist a finite uniform bound, i.e., any algorithm will potentially make arbitrarily large errors on any sufficiently large compact input domain.   
\end{remark}
How does the finding relate to our intuition and empirical evidence? Can we precisely identify input instances or input domains that lead to algorithmic failure? Regarding the former, we demonstrate, based on iterative processes with convergence guarantees, at what level the algorithmic breakdown potentially arises, which leads to concrete observations about the latter question. 

\subsection{Non-Effectiveness of Iterative Processes computing the Pseudoinverse}

It is important to stress that Theorem \ref{thm:non-approx} does not contradict the existence of procedures, which compute the pseudoinverse with provable convergence guarantees. 
For a non-zero matrix $A \in \C^{m\times n}$ with rank $p$, the iterative process $(A_k)_{k\in\N}$ based on
\begin{equation}\label{eq:ItPr}
    A_k :=
        \begin{cases}
            \alpha A^H &: \text{ if } k = 0,\\[3pt] 
            2 A_{k-1} - A_{k-1} A A_{k-1} &: \text{ if } k \in \N,
        \end{cases}
	\quad \text{ for some } 0 < \alpha < \frac{2}{\norm[2]{A}^2}, 
\end{equation}
converges to the pseudoinverse $A^\dagger$ via
\begin{equation}\label{eq:conv}
    \nnorm[2]{A^\dagger - A_k} \leq \frac{\norm[2]{A}}{\lambda_{p}(A^H A)} (1 - \alpha\lambda_{p}(A^H A))^{2^k} \quad \text{ for } k \in \N \cup \{0\},
\end{equation}
where $\norm[2]{\cdot}$ denotes the spectral norm and $\lambda_{p}(A^H A)>0$ the $p$-th largest eigenvalue of $A^H A$ \cite{BenIsrael66IterAlg}. In particular, the convergence rate obeys
\begin{equation}\label{eq:itCR}
	\nnorm[2]{A^\dagger - A_k} \leq \norm[2]{A} \|A^\dagger - A_{k-1}\|_2^2 \quad \text{ for } k \in \N.
\end{equation} 
This iterative process seemingly yields an effective algorithm to compute the pseudoinverse to any desired precision via the convergence of the sequence $(A_k)_{k\in\N}$. However, convergence is insufficient to guarantee an effective algorithm's existence. Indeed, the iteration step in \eqref{eq:ItPr} is expressed via the Borel-Turing computable function $f: \C^{n\times m} \times \C^{m\times n} \to \C^{n\times m}$ given by $f_B(X):=f(X, B) = 2 X - X B X$ as
\begin{equation}\label{eq:fIt}
    A_k = f_A(A_{k-1}) = f_A(f_A( \cdots f_A(A_0))) \quad \text{ for } k \in \N. 
\end{equation}
Hence, the sequence $(A_k)_{k\in\N}$ is a computable sequence provided that the initialization $A_0$ is a computable matrix since the composition of Borel-Turing computable functions is computable. Crucially, implementing the iterative process requires a stopping criterion that aborts the computation once the approximation is sufficiently close to $A^\dagger$. Ideally, the algorithm would take an error parameter $N \in \N$ as additional input and halt if, for some $m\in\N$, the approximation error of $A_m$ is smaller than $2^{-N}$. Therefore, an effective implementation takes a representation of any (non-zero) matrix $A\in\C^{m\times n}$ and an error parameter $N$ as input and outputs (an approximation of) $A^\dagger$ up to an error of $2^{-N}$.
The key point is that the initialization and the termination of the iteration are integral parts of the effective implementation. Otherwise, the iterative process could not be performed independently and effectively on digital hardware.

In the context of Turing machines the question of the effective implementation can be formalized in the following way: Does there exist a Borel-Turing computable function $G: \C^{m\times n} \to \C^{n\times m}$ such that $(A_k)_{k\in\N}$ computed according to \eqref{eq:fIt} with $A_0 = G(A)$ converges effectively to the pseudoinverse $A^\dagger$, i.e., does there exist a recursive function $e_G:\N \to \N$ such that, for any $A \in \C^{m\times n}$ with $G(A)=A_0$,
\begin{equation}\label{eq:Ginit}
    \nnorm[2]{A^\dagger - A_k} \leq \frac{1}{2^N} \quad \text{ holds true for all }  N \in \N \text{ and } k \geq e_G(N) ?
\end{equation}
Assuming that $G$ and $e_G$ exist entails that the representation of any non-zero matrix $A\in\C^{m\times n}$ can be effectively transformed into a representation of its pseudoinverse $A^\dagger$ via \eqref{eq:Ginit}, i.e., the function $f^\dag$ is Borel-Turing computable on $\C^{m\times n}\setminus\{0\}$. This contradicts Theorem \ref{thm:non-approx} stating that \eqref{eq:Ginit} can not even be realized for $N=1$ on a compact input set $\norm[F]{A}\leq \sqrt{2}$---our proof in the \nameref{sec:proofs} verifies that excluding the zero matrix from the input domain does not impact the result.

\begin{corollary}\label{cor:Effect}
    For $n,m\geq 2$, there does not exist a Borel-Turing computable function $G: \C^{m\times n} \to \C^{n\times m}$ to compute an initialization $A_0 = G(A)$ such that the sequence $(A_k)_{k\in\N}$, defined by
    \begin{equation*}
        A_k = 2 A_{k-1} - A_{k-1} A A_{k-1} \quad \text{ for } k \in \N,
    \end{equation*}
    converges effectively to the pseudoinverse $A^\dagger$.
\end{corollary}
Since the spectral norm is bounded by the (computable) Frobenius norm, fixing a computable function $G$, which provides a suitable initialization via \eqref{eq:ItPr}, is straightforward. However, the associated recursive function $e_G$ measuring the convergence rate via \eqref{eq:conv} does not exist because determining the rank of a matrix $A\in \C^{m\times n}$ and finding a (non-zero) lower bound on the $\rank(A)$-th largest eigenvalue of $A^H A$ can not be effectively computed without further information---hereby, determining the rank effectively is the prohibitive step. On the other hand, providing additional information about the input matrices by restricting the input domain, e.g., in the edge case to contain exactly one computable element, the iterative process yields an effective approach, as the next result shows.

\begin{proposition}\label{prop:comp_values}
   For any function I.-VI., i.e., $h \in \{f^\dag, f^{\nnorm{\dag}}, \kappa, f_{\text{lsq}}, f_\text{lsq-m},f_\text{lsq-n}\}$, and computable input $x \in \text{dom}(h)$, $h(x)$ is also computable.
\end{proposition}
\begin{remark}
        Corollary \ref{cor:Effect} and Proposition \ref{prop:comp_values} illustrate the contrast between effective initialization and effective convergence of the iterative process to obtain the pseudoinverse. If an effective initialization is implemented, then effective convergence is not feasible (Corollary \ref{cor:Effect}). Likewise, if effective convergence is guaranteed, then the corresponding initialization can not be performed effectively but needs to be manually crafted (Proposition \ref{prop:comp_values}).
\end{remark}
The previous discussion shows that any iterative process needs to rely on heuristics for initialization and stopping criteria. The heuristics may work reasonably well for specific classes of matrices, but they can not guarantee correctness (in the effective sense) in the general case. In particular, the proof of Proposition \ref{prop:comp_values} is not constructive and can therefore not be generalized to an algorithm taking arbitrary computable inputs since specific properties of the given instance need to be provided to the associated algorithm. However, is it possible to at least establish effective approaches (and not only heuristics) to tackle the problems I.-VI. for specific classes of matrices? This is indeed feasible as we establish next.

\subsection{Full description of computability breakdown}
Can we identify the properties of matrices that lead to the computability breakdown and thereby identify domains that allow for effective treatment? The proof of Theorem \ref{thm:non-approx} and the short analysis of the SVD approach in Section \ref{sec:BMnonComp} indicate that non-computability indeed arises in specific circumstances. The following result and its implication show that non-computability can be prevented in certain instances, although it can not be circumvented for sufficiently general domains.  
\begin{theorem}\label{thm:class}
   For any function I.-VI., i.e., $h \in \{f^\dag, f^{\nnorm{\dag}}, \kappa, f_{\text{lsq}}, f_\text{lsq-m},f_\text{lsq-n}\}$, there exists a Turing machine that transforms a representation of any computable input $x \in \text{dom}(h)$ into a representation for $h(x)$ provided that each member of the representation of $x$ has the same rank as $x$ itself.
\end{theorem}
\begin{remark}
    Note that Theorem \ref{thm:class} does not imply Borel-Turing computability. Indeed, the non-compliance with the rank property is exactly the condition exploited in the proof of the non-approximability result in Theorem \ref{thm:non-approx}.    
\end{remark}
Theorem \ref{thm:non-approx} together with Theorem \ref{thm:class} fully characterize the algorithmic approximability of the pseudoinverse and the related functions II.-VI.: 
\begin{itemize}
    \item If the rank property is satisfied, then the rank of any input matrix can be directly computed via its given representation. Subsequently, the pseudoinverse can be computed to any desired accuracy by the constructive algorithm described in the proof of Theorem \ref{thm:class}. 
    \item In certain settings, the rank property can be effectively guaranteed, i.e., there exists an effective algorithm that transforms any given representation into a representation meeting the rank property, e.g., the set of full rank matrices in $\C^{m\times n}$ or the set $\C^{m\times n}_\varepsilon :=\{A\in \C^{m\times n} : |A_{ij}| > \varepsilon \text{ or } A_{ij}=0\,\, \forall i,j\}$ for any $\varepsilon>0$ (see \nameref{sec:proofs}, in particular, Theorem \ref{thm:FullRank} and Proposition \ref{prop:EpsSet} for details). 
    Therefore, the pseudoinverse mapping is Borel-Turing computable on the respective domains. 
    However, in both cases, the membership to the given set can not be effectively decided, i.e., for general input domains the best one can hope for is a pre-processing step that verifies if a given input belongs to a 'benign' class of matrices (but does not effectively distinguish between benign and non-benign classes), i.e., only semi-decides the membership property.
    \item The algorithmic non-approximability arises for input domains being arbitrarily close to or containing matrices with distinct rank (if no further conditions are imposed on the input instances). This is, in some sense, an expected behaviour because it coincides with the discontinuities of the pseudoinverse mapping \cite{Rakočević1997}.
    In these cases, the construction in the proof of Theorem \ref{thm:FullRank} can be encoded in the considered input domain, leading to a non-approximability result. Hence, the key feature of the set of full rank matrices and $\C^{m\times n}_\varepsilon$ is preventing a (smooth) rank transition by design---the set of full rank matrices is open, whereas $\C^{m\times n}_\varepsilon$ artificially imposes via $\varepsilon$ a testable condition on rank transitioning. In contrast, the set of matrices with fixed rank smaller than $\min(m,n)$ also contains only matrices with equal rank, however, the set is closed so that for instances on the boundary, there exist representations solely consisting of matrices with a distinct rank.
\end{itemize}

\section*{Acknowledgement}
Holger Boche acknowledges the financial support by the BMBF in the programme of ”Souverän. Digital. Vernetzt”, research HUB 6G-life, project identification number: 16KISK002, and the financial support by the BMBF Quantum Projects QUIET, Grant 16KISQ093, QD-CamNetz, Grant 16KISQ077, and QuaPhySI, Grant 16KIS1598K. Furthermore, he acknowledges funding by the German Research Foundation (DFG, Deutsche Forschungsgemeinschaft) as part of Germany’s Excellence Strategy – EXC 2050/2 – Project ID 390696704 – Cluster of Excellence “Centre for Tactile Internet with Human-in-the-Loop” (CeTI) of Technische Universität Dresden. 
H. Boche was also partially supported by the project “Next Generation AI Computing (gAIn)”, funded by the Bavarian Ministry of Science and the Arts and the Saxon Ministry for Science, Culture, and Tourism.

This work of G. Kutyniok was supported in part by the Konrad Zuse School of Excellence in Reliable AI (DAAD), the Munich Center for Machine Learning (MCML) as well as the German Research Foundation under Grants DFG-SPP-2298, KU 1446/31-1 and KU 1446/32-1. Furthermore, G. Kutyniok acknowledges additional support by the project ”Next Generation AI Computing (gAIn)”, funded by the Bavarian Ministry of Science and the Arts and the Saxon Ministry for Science, Culture, and Tourism. % as well as by the Hightech Agenda Bavaria.  

\printbibliography

\section*{Appendix}
\label{sec:proofs}

\subsection*{Proof of Theorem \ref{thm:non-approx}} 

The proof of Theorem \ref{thm:non-approx} relies on general non-approximability conditions that we present in the following lemma. They were introduced in \cite{Boche2023InvProbDig}, and we provide a slightly adapted version for clarity, yet we refer to \cite{Boche2023InvProbDig} for a proof, which remains structurally unchanged.
\begin{lemma}\label{lemma:non-approx}
    Let $n,m\in\N$ and consider the map $f : X \to Y$, where $X \subset \C^{m\times n}$ or $X \subset \C^{m\times n}\times \C^m$ and $Y \subset \C^n$, $Y \subset \C^{n\times m}$ or $Y \subset \R$. Denote by $\norm[X]{\cdot}$ a norm on $X$ and by $\norm[Y]{\cdot}$ a norm on $Y$, which is Banach-Mazur computable (as a function). Further, suppose that there exists a computable sequence $(x_{n})_{n\in \N}\subset X$ satisfying the following conditions:
    \begin{enumerate}[label=(\alph*)]
        \item There exists $x^\ast \in X$ such that $\norm[X]{x_{n} - x^\ast} \leq 2^{-n}$ for all $n \in \N$. 
        \item There exists $\eta > 0, \eta \in \Q$ such that $\inf_{n \in \N} \norm[Y]{f(x_n) - f(x^\ast)} > \eta$.
    \end{enumerate}
    Then any function $\hat{f}: X \to Y$ with
    \begin{equation*}
        \sup_{x \in X} \norm[Y]{f(x) - \hat{f}(x)} < \frac{\eta}{4}
    \end{equation*}
    is not Banach-Mazur computable.   
\end{lemma}
By employing Lemma \ref{lemma:non-approx}, we prove Theorem \ref{thm:non-approx} via a suitable construction.
\begin{proof}[Proof of Theorem \ref{thm:non-approx}]
    We show that the conditions of Lemma \ref{lemma:non-approx} are satisfied for each function I.-VI. and thereby derive the given non-approximability bounds in Banach-Mazur sense. Consider for fixed dimensions $m\geq 2$, $n=2$ and $\varepsilon\geq 0$ the matrix-vector pairs 
    \begin{equation*}
        A_\varepsilon := \begin{pmatrix} 1 & 0 \\ 0  & \varepsilon \\ 0 & 0 \\ \vdots & \vdots \\ 0 & 0 \end{pmatrix} \in \C^{m \times n} \quad \text{ and } \quad b := \begin{pmatrix} 1 \\ 1 \\ 0 \\ \vdots \\ 0\end{pmatrix} \in \C^m,
    \end{equation*}
    where $m-2$ zero rows are appended, respectively. The general case $n > 2$ follows via the same proof strategy by appending $n-2$ zero-columns to $A_\varepsilon$. Hence, without loss of generality, we only describe the case $n=2$. We start with the following observations:
    \begin{itemize}
        \item For $\varepsilon>0$, the columns of $A_\varepsilon$ are linearly independent, i.e., $A^T_\varepsilon A_\varepsilon$ is invertible, so that 
        \begin{equation}\label{eq:Aepsdag}
            A^\dagger_\varepsilon = (A^T_\varepsilon A_\varepsilon)^{-1} A^T_\varepsilon = \begin{pmatrix} 1 & 0 & 0 & \dots & 0 \\ 0 & \varepsilon^{-1} & 0 & \dots & 0  \end{pmatrix} \text{  }  \text{  and  } \text{  }
            \hat{x}_{(A_\varepsilon,b)} := A^\dagger_\varepsilon b = 
            \begin{pmatrix}
                1 \\ \varepsilon^{-1}
            \end{pmatrix}.
        \end{equation}
        \item Similarly for $\varepsilon=0$, we obtain
        \begin{equation}\label{eq:A0dag}
            A^\dagger_0 = \begin{pmatrix} 1 & 0 & 0 & \dots & 0 \\ 0 & 0 & 0 & \dots & 0  \end{pmatrix} \quad \text{ and } \quad \hat{x}_{(A_0,b)} := A^\dagger_0 b = \begin{pmatrix}
                1 \\ 0
            \end{pmatrix}.
        \end{equation}
        \item Finally, notice that 
        \begin{equation}\label{eq:AdiffNorm}
            \norm[F]{A_\varepsilon - A_0} = \varepsilon \quad \text{ and } \quad  \nnorm[F]{A^\dagger_\varepsilon - A^\dagger_0} = \varepsilon^{-1} \text{ for } \varepsilon >0.
        \end{equation}
    \end{itemize}
    For the sequence $(\varepsilon_n)_{n\in\N}$ given by $\varepsilon_n = 2^{-n}$, it is immediate to verify that the associated sequence of matrices $(A_{2^{-n}})_{n\in\N}$ is computable. Additionally, by applying \eqref{eq:AdiffNorm}, we obtain for $n\in\N$ that
    \begin{equation*}
        \norm[F]{A_{2^{-n}} - A_0} = 2^{-n} \quad \text{and} \quad \norm[F]{f^\dag(A_{2^{-n}}) - f^\dag(A_0)} = \nnorm[F]{A^\dagger_{2^{-n}} - A^\dagger_0} = 2^{n}.   
    \end{equation*}
    Thus, for any $c>0$ there exists $k_c\in\N$ such that
    \begin{equation*}
        \norm[F]{f^\dag(A_{2^{-(k_c+\ell)}}) - f^\dag(A_0)} > 4c \quad \text{ for all } \ell \in \N.
    \end{equation*}
    Since $\norm[F]{\cdot}$ is Banach-Mazur computable, the computable sequence $(A_{2^{-(n+k_c)}})_{n\in\N}$ satisfies the conditions in Lemma \ref{lemma:non-approx} with $\eta = 4c$. Since this holds for any $c>0$, $\eta$ can be chosen arbitrarily large. 
    Furthermore, it suffices to consider the compact input domain $A_{\text{compact}}\:= \{ A \in \C^{m\times n} : \norm[F]{A} \leq \sqrt{2} \} \subset \C^{m\times n}$ and invoke Lemma \ref{lemma:non-approx} on $A_{\text{compact}}$ since 
    \begin{equation*}
        \norm[F]{A_0}, \sup_{n\in\N, k_c} \norm[F]{A_{2^{-(n+k_c)}}}\leq \sqrt{2}, \text{ i.e., } \{A_0\} \cup \{A_{2^{-(n+k_c)}}\}_{n\in\N} \subset A_{\text{compact}}\quad  \forall k_c.
    \end{equation*}
    For the treatment of $f^{\nnorm{\dag}}$, a similar reasoning together with the Banach-Mazur computability of $\abs{\cdot}$ and the observation
    \begin{equation*}
        \abs{f^{\nnorm{\dag}}(A_{2^{-n}}) - f^{\nnorm{\dag}}(A_0)} = \abs{\nnorm[F]{A^\dagger_{2^{-n}}} - \nnorm[F]{A^\dagger_0}} = \sqrt{1^2 + (2^{n})^2} - 1 \geq n \quad \text{ for } n\in\N
    \end{equation*}
    suffice.
    Regarding $\kappa$, applying \eqref{eq:Aepsdag} and \eqref{eq:A0dag}, yields that 
    \begin{equation}\label{eq:CondF}
        \norm[F]{A_\varepsilon} \nnorm[F]{A^\dagger_\varepsilon} = \sqrt{1 + \varepsilon^2} \sqrt{1+\varepsilon^{-2}} \text{ for } \varepsilon >0  \quad \text{ and } \quad \norm[F]{A_0} \nnorm[F]{A^\dagger_0} = 1.     
    \end{equation}
    Consequently, we can conclude that for $n\in \N$
    \begin{align*}
        \abs{\kappa(A_{2^{-n}}) - \kappa(A_0)} &= \abs{\norm[F]{A_{2^{-n}} } \nnorm[F]{A^\dagger_{2^{-n}} } - \norm[F]{A_0} \nnorm[F]{A^\dagger_0} } = \sqrt{1 + 2^{-2n}} \sqrt{1+2^{2n}} - 1 \\
        &\geq  \sqrt{1+2^{2n}} - 1 \geq n.
    \end{align*}
    Therefore, analogously to the previous constructions, the claim with respect to $\kappa$ follows from Lemma  \ref{lemma:non-approx}.
    Finally, we prove that the applied norm---the Frobenius norm---in the definitions of $f^{\nnorm{\dag}}$ and $\kappa$ may be replaced with any other norm. We explicitly show the statement for $\kappa$, however, the other case $f^{\nnorm{\dag}}$  follows analogously. 
    Let $\norm[c]{\cdot}$ be an arbitrary matrix norm and denote by $\tilde{\kappa}: \C^{m\times n} \to \R$ the mapping $\tilde{\kappa}(A)  := \norm[c]{A} \nnorm[c]{A^\dagger}$.
    Since all norms are equivalent on a finite-dimensional space, there exist constants $c,C,d,D > 0$ (w.l.o.g. $C,D\geq 1$) such that, for any matrix $M \in \C^{m\times n}$, we have
    \begin{equation*}
        c \norm[F]{M} \leq \norm[c]{M} \leq C \norm[F]{M} \quad  \text { and } \quad d \nnorm[F]{M^\dagger} \leq \nnorm[c]{M^\dagger} \leq D \nnorm[F]{M^\dagger}. 	
    \end{equation*}
    Thus, applying \eqref{eq:CondF} yields
    \begin{align}\label{eq:cnorm}
        \norm[c]{A_\varepsilon} \norm[c]{A^\dagger_\varepsilon}  &\geq c d \norm[F]{A_\varepsilon} \norm[F]{A^\dagger_\varepsilon} = c d \sqrt{1 + \varepsilon^2} \sqrt{1+\varepsilon^{-2}} \text{ for } \varepsilon >0 \quad \text{ and }  \nonumber \\
        \quad \norm[c]{A_0} \nnorm[c]{A^\dagger_0} &\leq C D \norm[F]{A_0} \nnorm[F]{A^\dagger_0} = C D.
    \end{align}
    Choose $r \in \N$ such that $c d 2^{r} > C D$ and observe that for $n\in\N$ we obtain  via \eqref{eq:cnorm} 
    \begin{align*}
        \abs{\tilde{\kappa}(A_{2^{-(r+n)}}) - \tilde{\kappa}(A_0)} &= \abs{\nnorm[c]{A_{2^{-(r+n)}} } \nnorm[c]{A^\dagger_{2^{-(r+n)}} } - \norm[c]{A_0} \nnorm[c]{A^\dagger_0} } \nonumber \\
        &\geq \abs{\nnorm[c]{A_{2^{-(r+n)}} } \nnorm[c]{A^\dagger_{2^{-(r+n)}} }} - \abs{\norm[c]{A_0} \nnorm[c]{A^\dagger_0} }\nonumber \\
        &\geq c d \sqrt{1 + 2^{-2(r+n)}} \sqrt{1+2^{2(r+n)}}  - C D \\
        &\geq cd 2^{r+n} - C D \geq C D (2^n -1) \geq 2^{n-1}
    \end{align*}
    Hence, invoking once again Lemma \ref{lemma:non-approx} with the sequence $(A_{2^{-(r+n)}})_{n \in \N}$ and the same considerations as before with respect to the input domain gives the claim for arbitrary computable norm, which finishes the proof of $(i)$.
    
    For $(ii)$, define the (Banach-Mazur computable) norm $\norm[(F,2)]{(\cdot, \cdot) }: \C^{m\times n} \times \C^m \to \R$ by $\norm[(F,2)]{(A,b)} := \norm[F]{A} + \norm[2]{b}$ and note that \eqref{eq:AdiffNorm} implies 
    \begin{equation}\label{eq:Ab}
        \norm[(F,2)]{(A_\varepsilon,b) - (A_0,b)} = \norm[F]{A_\varepsilon - A_0} + \norm[2]{b - b}  = \varepsilon.
    \end{equation}
    Moreover, we conclude via \eqref{eq:Aepsdag} and \eqref{eq:A0dag} that
    \begin{equation}\label{eq:xAb}
        \norm[2]{A_\varepsilon \hat{x}_{(A_\varepsilon,b)} - b} = 0 \text{ for } \varepsilon >0 \quad \text{ as well as } \quad \norm[2]{A_0 \hat{x}_{(A_0,b)} - b}  =  1.
    \end{equation}
    Applying the computable sequence of matrix-vector pairs $((A_{2^{-n}},b))_{n\in\N}$, we can invoke Lemma  \ref{lemma:non-approx} with any $0 <\eta <1$ since \eqref{eq:Ab} yields that
    \begin{equation*}
        \norm[(F,2)]{(A_{2^{-n}},b) - (A_0,b)}= 2^{-n}
    \end{equation*}
    and \eqref{eq:xAb} gives that
    \begin{equation*}
        \abs{f_{\text{lsq}}(A_{2^{-n}},b) - f_{\text{lsq}}(A_0,b)} = \abs{\nnorm[2]{A_{2^{-n}} \hat{x}_{(A_{2^{-n}},b)} - b} - \nnorm[2]{A_0 \hat{x}_{(A_0,b)} - b} } = 1.
    \end{equation*}
    Finally, observing that $\norm[2]{b} \leq \sqrt{2}$ implies that it suffices to consider the domain $A_{\text{compact}} \times \{ v \in \C^m : \norm[2]{v} \leq \sqrt{2} \} \subset \C^{m\times n} \times \C^m$, which shows $(ii)$.

    Regarding $(iii)$, this follows analogously from the previous constructions together with the observations derived via \eqref{eq:Aepsdag} and \eqref{eq:A0dag} that for any $n\in\N$ 
    \begin{equation*}
        \norm[2]{f_\text{lsq-m}(A_{2^{-n}},b) - f_\text{lsq-m}(A_0,b)} = \nnorm[2]{\hat{x}_{(A_{2^{-n}},b)} - \hat{x}_{(A_0,b)}} = 2^{n} 
    \end{equation*}
    and
    \begin{equation*}
        \abs{f_\text{lsq-n}(A_{2^{-n}},b)-\Psi_{\text{norm}}(A_0,b)} = \abs{ \nnorm[2]{\hat{x}_{(A_{2^{-n}},b)}} - \nnorm[2]{\hat{x}_{(A_0,b)}} } = \sqrt{1 + 2^{2n}} - 1 \geq n.  
    \end{equation*}     
\end{proof}

\subsection*{Proof of Proposition \ref{prop:comp_values}}

Here, we provide the proof of Proposition \ref{prop:comp_values}.
\begin{proof}[Proof of Proposition \ref{prop:comp_values}]
    For the zero matrix as input, the statement trivially holds. Moreover, Remark \ref{rm:elementaryfct} implies that if the pseudoinverse is a computable matrix for a given input $x \in \text{dom}(h)$, then $h(x)$ is a computable real/vector/matrix as well. Hence, it is left to show that for any non-zero, computable matrix $A \in \C_c^{m\times n}$ the corresponding pseudoinverse $A^\dagger$ is a computable matrix. To that end, initialize $A_0= \beta A^H$ with $\beta \in \R_c$ satisfying the condition in \eqref{eq:ItPr}---$\beta$ indeed exists since $\R_c$ is dense in $\R$---so that $(A_k)_{k\in\N}$ is a computable sequence converging to $A^\dagger$. Thus, there exists $k^\ast \in \N$ such that
    \begin{equation*}
        \nnorm[2]{A^\dagger - A_{k^\ast}} \leq 
            \begin{cases}
                \frac{1}{2} &: \text{ if } \nnorm[2]{A} \leq 1,\\[3pt] 
                \frac{1}{2\norm[2]{A}} &: \text{ if } \nnorm[2]{A} > 1.
            \end{cases}
    \end{equation*}
    It follows that the computable sequence $(A_{k^\ast + k})_{k\in\N}$ satisfies for any $k\in\N$ that
    \begin{equation*}
       \nnorm[2]{A^\dagger - A_{k^\ast + k}} \leq \norm[2]{A}^{2^k-1} \|A^\dagger - A_{k^\ast}\|_2^{2^k}  \leq \frac{1}{2^k},
    \end{equation*}
    where we have applied \eqref{eq:itCR}. 
    Hence, the computable sequence $(A_{k^\ast + k})_{k\in\N}$ converges effectively to $A^\dagger$ and therefore $A^\dagger$ is a computable matrix.
\end{proof}

\subsection*{Proof of Theorem \ref{thm:class}}\label{app:thmClass}

We start by analyzing a specific input domain, the set of full rank matrices, and show that an effective approach to compute the pseudoinverse can be constructed in this setting. We continue by introducing a characterization of input domains that enable the existence of effective algorithms and is also the basis for the proof of Theorem  \ref{thm:class}. Lastly, we exemplarily demonstrate the application of the characterization to derive an effective algorithm for a specific input set.

\subsubsection*{Full Rank Matrices}
We consider invertible matrices $A\in \text{GL}_n(\C)$ first. For any $A\in \text{GL}_n(\C)$, we have $A^\dag= A^{-1}$ and there exists an effective approach to compute the inverse provided that the input instances are restricted to invertible matrices. This can be straightforwardly verified via Gaussian elimination, however, we present a different approach, which will also be useful for more general input domains below. 
\begin{lemma}\label{lm:Inv}
    For any $n\in\N$, $ \Psi\colon \text{GL}_n(\C)\to \text{GL}_n(\C),\, A\mapsto A^{-1}$ is Borel-Turing computable.
\end{lemma}
\begin{proof}
    For any invertible matrix $A\in \C^{n\times n}$, i.e., $\det(A)\neq 0$, the inverse can be represented as
    \begin{equation}\label{eq:InvFormula}
        A^{-1}=\frac{-1}{\det(A)}(\alpha_0I+\alpha_1A+...+\alpha_nA^{n-1}),
    \end{equation}
    where $\alpha_i$ are the coefficients of the characteristic polynomial $\mathcal{X}_A$ of A, i.e.
    \begin{equation*}
        \mathcal{X}_A(t)=\det(t I - A) = \alpha_n t^n+ \cdots + \alpha_1t+ \alpha_0.
    \end{equation*}
    As described in for instance \cite{gantmacher1980theory}, the coefficients $(\alpha_0, \dots, \alpha_n)$ can be computed via
    \begin{align}\label{eq:coeffComp}
        \alpha_n&=1 \nonumber \\
        \alpha_{n-k} &= \frac{(-1)^{n-k}}{k!} B_k(s_1,\dots,s_k)\quad  \text{ for } k=1,\dots, n-1 \nonumber\\
         \alpha_0 &= \frac{1}{n!} B_n(s_1,\dots,s_n)=(-1)^n \det(A),    
    \end{align}
    where $s_i=-(i-1)\tr(A^i)$ for $i=1,\dots,n$ and $B_\ell$ for $\ell=1,\dots,n$ is the $\ell$-th complete exponential Bell polynomial \cite{Bell34Poly} 
    \begin{equation*}
        B_\ell(x_1,\dots,x_\ell)=\ell! \sum_{\substack{j_m\in \Z, j_m\geq 0 : \\ \sum_{m=1}^\ell m j_m=\ell}}  \prod_{i=1}^\ell \frac{x_i^{j_i}}{(i!)^{j_i}j_i!}.
    \end{equation*}
    It follows that $B_\ell: \C^\ell \to \C$ considered as a map is Borel-Turing computable since $B_\ell$ is characterized by computable operations with the condition under the sum effectively verifiable for all possible choices of $(j_1, \dots, j_\ell)\subset \{0,\dots,\ell\}^\ell$. Therefore, observing that the map $A\mapsto (s_1,\dots,s_n)$ is Borel-Turing computable together with \eqref{eq:coeffComp} immediately implies that the determinant function $\det: \C^{n\times n} \to \C$ and the function $g: \C^{n\times n} \to \C^{n\times n}$ given by $B\mapsto \beta_0I+\beta_1B+\dots+\beta_nB^{n-1}$, where $\beta_1,\dots, \beta_n$ are the coefficients of the characteristic polynomial of $B$, are Borel-Turing computable. 
    
    Note that for invertible matrices $A$, the inverse is given by $\det(A)^{-1} g(A)$ according to \eqref{eq:InvFormula}. Hence, an effective algorithm transforming each representation of any invertible matrix into a representation of its inverse exists, provided that the input representation solely consists of invertible matrices. The condition is necessary since \eqref{eq:InvFormula} is only valid for invertible matrices (otherwise, the formula would include a division by zero). However, showing Borel-Turing computability of the inverse mapping $\Psi$ requires the described transformation to hold for any representation, including ones with non-invertible members.
    
    We can solve this issue by introducing an effective processing step that takes an arbitrary representation and transforms it into a representation consisting of invertible matrices before tackling the inversion. Informally, by iteratively increasing the accuracy of a representation $(A_k)_{k\in\N}$ of $A \in \text{GL}_n(\C_c)$, we can effectively find $m^\ast \in\N$ such that both $|\det(A) - \det(A_{m^\ast})| < 2^{-(m^\ast+1)}$ and $|\det(A_{m^\ast})| > 2^{-m^\ast}$ hold, which implies that $|\det(A_{m^\ast + m})| > 0$ for all $m\in\N$, i.e., $(A_{k+m^\ast})_{k\in\N}$ is a representation of $A$ consisting of invertible matrices. The process is guaranteed to succeed since $\det$ is Borel-Turing computable on $\C^{n\times n}$ and $\text{GL}_n(\C)$ is open in $\C^{n\times n}$. Therefore, we can conclude that $\Psi$ is Borel-Turing computable. 
\end{proof}

The invertible case is relevant because computing the pseudoinverse in the more general case of matrices $A\in \C^{m\times n}$ with full rank boils down to inversion \cite{BenIsrael2003GenInv}: 
\begin{itemize}
    \item either $A^\dag=(A^H A)^{-1}A^H$ if $A$ has linearly independent columns,
    \item or $A^\dag=A^H(AA^H)^{-1}$ if $A$ has linearly independent rows.
\end{itemize}

\begin{theorem}\label{thm:FullRank}
    For any $n,m\in\N$, $ \Psi\colon \text{FR}_{m,n}(\C)\to \text{FR}_{n,m}(\C),\, A\mapsto A^\dagger$ is Borel-Turing computable with $\text{FR}_{m,n} := \{ A\in\C^{m\times n} :\,  \rank(A)=\min\{m,n\} \}$.
\end{theorem}
\begin{proof}
    Note that for $m=n$, we are back in the setting of the above lemma. The remaining possibilities $m<n$ and $m>n$ follow along the same lines, so we will only present the details for $m>n$.
    In this case, the columns of the considered matrices are linearly independent so that the pseudoinverse of $A\in \text{FR}_{m,n}$ is given by $A^\dag=(A^H A)^{-1}A^H$, which via Lemma \ref{lm:Inv} we identify as a composition/product of computable operations. Hence, an effective algorithm transforming each representation of any matrix with full rank into a representation of its pseudoinverse exists, provided that the input representation consists of matrices with full rank. To attain Borel-Turing computability of $\Psi$, we introduce an effective process that takes an arbitrary representation and transforms it into a representation consisting of full rank matrices. 

    This is achieved with the same technique as in the previous lemma applied to the matrices $(A_k^H A_k)_{k\in\N}$ for a given representation $(A_k)_{k\in\N}$, taking advantage of the fact that $A^H A \in \text{GL}_n(\C)$ by the full rank condition. Hence, $m^\ast\in\N$ can be effectively computed so that $(A_{k+m^\ast})_{k\in\N}$ is a representation of $A$ consisting of full rank matrices, and therefore, the claim follows.       
\end{proof}

\subsubsection*{Representations with Fixed Rank}
Enlarging the input domain beyond matrices with full rank while maintaining the Borel-Turing computability of the mapping on the pseudoinverse requires rather specific classes of matrices for two reasons. If matrices with different ranks are admissible, then, in general, the construction in the proof of Theorem \ref{thm:non-approx} allows us to invoke Lemma \ref{lemma:non-approx}, leading to the algorithmic non-approximability of the pseudoinverse mapping on the given domain. Restricting the domain to matrices of fixed rank (except for the previously treated full rank case) also leads to algorithmic non-approximability due to input instances on the domain boundary---the set of matrices with fixed rank is closed. 
A different approach is to relax the computability requirement and thereby extend the relaxed notion to more general input domains.

\begin{theorem}\label{thm:pos}
    There exists a Turing machine that transforms a representation $(A_k)_{k\in\N}$ of any matrix $A\in\C_c^{m\times n}$ into a representation of $A^\dag$ provided that $\rank(A_k) =\rank(A)$ for all $k\in\N$. 
\end{theorem}
\begin{proof}
For any $A\in \C^{m\times n}$ with $\rank(A) = p < \min\{m,n\}$, the pseudoinverse satisfies $A^\dag=B A^H$ \cite{Penrose55GenInv} with 
\begin{equation}\label{eq:ConstructB}
    B=-\alpha_r^{-1}\big(\alpha_{r+1}I+...+\alpha_n(A^HA)^{n-r-1}\big),
\end{equation}
where $r$ is the smallest index of a non-zero coefficient $\alpha_i$ given by
\begin{equation*}
    \mathcal{X}_{A^HA}(A^HA) = \alpha_0 I + \alpha_1 A^HA +\dots + \alpha_n (A^HA)^n = \alpha_1 A^HA + \dots +\alpha_n (A^HA)^n.    
\end{equation*}
The index $r$ exists since $A^HA,(A^HA)^2,\dots,(A^H A)^n$ are linearly dependent and we used that $\alpha_0 =(-1)^n \det(A^H A)$ via \eqref{eq:coeffComp} combined with the fact that $\det(A^H A)=0$ since $\rank(A^H A) = \rank(A)=p < n$.
Since $A^HA$ is Hermitian, it is diagonalizable and has $j$ (distinct) eigenvalues $\lambda_1,\dots,\lambda_j \in \R$ with algebraic multiplicities $a_1 +\dots+ a_j=n$.
We know that $A^HA$ is not invertible, so $\lambda_1=0$ is an eigenvalue with algebraic multiplicity (equaling geometric multiplicity) 
\begin{equation*}
    a_1= \dim \ker(A^H A) = n - \rank(A^HA) = n - \rank(A)= n - p.    
\end{equation*}
Therefore, we can write
\begin{align*}
    \mathcal{X}_{A^HA}(t) &=\prod\limits_{i=1}^j (t - \lambda_i)^{a_i}= t^{n - p}\prod\limits_{i=2}^j(t - \lambda_i)^{a_i}  \\
    &=t^{n-p}\prod\limits_{i=2}^j (-\lambda_i) ^{a_i} + t^{n-p+1}\alpha_{n-p+1} +\dots + t^n \alpha_n, 
\end{align*}
i.e., the first non-zero coefficient $\alpha_r=\prod_{i=2}^j(-\lambda_i)^{a_i}\neq 0$ of $\mathcal{X}_{A^HA}$ is attained for $r=n-p$. Subsequently,  the coefficients $\alpha_r, \dots, \alpha_n$ can be effectively computed via the formula \eqref{eq:coeffComp}. Therefore, given a representation $(A_k)_{k\in\N}$ of $A$ with $\rank(A_k) = \rank(A)=p$, the sequence $(A_k^\dag)_{k\in\N} = (B_k A_k^H)_{k\in\N}$ converges effectively to $B A^H = A^\dag$, where each $B_k \in \C^{n\times n}$ and $B$ are given according to \eqref{eq:ConstructB} with initial matrices $A_k$ and $A$ corresponding to coefficients $\alpha^k=(\alpha^k_j)_{j=r}^n$ and $\alpha=(\alpha_j)_{j=r}^n$, respectively.
This follows since the structure of the formula \eqref{eq:ConstructB} remains unchanged for each $B_i$ and $B$ due to the rank condition; in particular, the coefficients with the same indices $r, \dots, n$ appear. Finally, observing that \eqref{eq:coeffComp} implies the effective convergence of the coefficient sequence $(\alpha^k)_{k\in\N}$ to $\alpha$ and consequently of $(B_k)_{k\in\N}$ to $B$, gives the effective convergence of $(A_k^\dag)_{k\in\N} = (B_k A_k^H)_{k\in\N}$ to $B A^H = A^\dag$.

Hence, assuming that one can effectively determine the rank of an arbitrary matrix $M\in \C_c{m\times n}$ given a representation $(M_k)_{k\in\N}$, the claim follows since $M^\dag$ can then be effectively computed by above algorithm ($\rank(M)< \min\{m,n\}$) or by invoking Theorem  \ref{thm:FullRank} ($\rank(M)= \min\{m,n\}$). Indeed, for representations satisfying the condition $\rank(M)=\rank(M_k)$ for all $k\in\N$, it suffices to compute the rank of any $M_k$. This can be achieved by observing that the rank is the largest order $q \leq \min\{m,n\}$ of any non-zero minor of $M_k$, i.e., the largest $q$ so that a matrix obtained from $M_k$ by removing $m-q$ rows and $n-k$ columns has a non-zero determinant. Since there are only finitely many minors of a given matrix, it suffices to compute iteratively the determinant of each minor and compare to zero, which can be effectively implemented because of the Borel-Turing computability of the determinant function (see the proof of Lemma \ref{lm:Inv}) and and the feasibility of effective comparisons with zero for rational numbers---the determinant of any minor of $M_k$ is rational since $M_k \in \Q^{m\times n}$. 
\end{proof}

Theorem \ref{thm:class} is now a direct consequence of our previous considerations. 

\begin{proof}[Proof of Theorem\ref{thm:class}]
    Proposition \ref{thm:pos} shows that the claim is true for $h=f^\dag$. The remaining cases $h \in \{f^{\nnorm{\dag}}, \kappa, f_{\text{lsq}}, f_\text{lsq-m},f_\text{lsq-n}\}$ follow directly by observing that $h$ can be expressed via a composition of $f^\dag$ and computable functions (via Remark \ref{rm:elementaryfct}).     
\end{proof}

By verifying that the condition in Theorem \ref{thm:pos} is satisfied (or can be effectively guaranteed) for every representation of input instances, one can establish Borel-Turing computability on the given input domain---this approach was already the basis of proving Lemma \ref{lm:Inv} and Theorem \ref{thm:FullRank}; now we apply it exemplarily to an input domain without rank restrictions. 

\begin{proposition}\label{prop:EpsSet}
    For any $n,m\in\N$ and $\varepsilon>0$, $ \Psi\colon \C^{m\times n}_\varepsilon \to \C^{n\times m},\, A\mapsto A^\dagger$, is Borel-Turing computable with $\C^{m\times n}_\varepsilon :=\{A\in \C^{m\times n} : |A_{ij}| > \varepsilon \text{ or } A_{ij}=0\,\, \forall i,j\}$.
\end{proposition}
\begin{proof}
    By definition, the zero and non-zero elements of any $A\in\C^{m\times n}_\varepsilon$ can be effectively identified via any representation $(A_k)_{k\in\N}$ when considering a sufficiently high accuracy (e.g. $\tfrac{\varepsilon}{2}$), i.e., a sufficiently large index $k^\ast$. Subsequently, setting 
    \begin{equation*}
        (A_k^\prime)_{ij} =  \begin{cases}
            (A_{k+k^\ast})_{ij}, &\text{ if } A_{ij} \neq 0, \\
            0 , &\text{ if } A_{ij} = 0
        \end{cases},    
    \end{equation*}
    yields a representation $(A^\prime_k)_{k\in\N}$ of $A$ with $\rank(A^\prime_k) = \rank(A)$ so that Theorem \ref{thm:pos} gives the result.
        
\end{proof}

%\printbibliography
\end{document}